\documentclass{article}
\usepackage{amssymb,amsmath}
\usepackage{graphicx}

\def\qed{\hfill {\hbox{${\vcenter{\vbox{               %HOLLOW SQUARE
   \hrule height 0.4pt\hbox{\vrule width 0.4pt height 6pt
   \kern5pt\vrule width 0.4pt}\hrule height 0.4pt}}}$}}}

\def\tr{\triangleright}

\newtheorem{theorem}{Theorem}
\newtheorem{definition}{Definition}
\newtheorem{lemma}[theorem]{Lemma}
\newtheorem{proposition}[theorem]{Proposition}
\newtheorem{corollary}[theorem]{Corollary}
\newtheorem{example}{Example}

\newtheorem{remark}{Remark}

\newenvironment{proof}[1][Proof]{\smallskip\noindent{\bf #1.}\quad}%
{\qed\par\medskip}

\date{}

\title{\Large \textbf{Link invariants from finite Coxeter racks}}

\author{Sam Nelson \and Ryan Wieghard}

\begin{document}
\maketitle

\begin{abstract}
We study Coxeter racks over $\mathbb{Z}_n$ and the knot
and link invariants they define. We exploit the module
structure of these racks to enhance the rack counting invariants
and give examples showing that these enhanced invariants
are stronger than the unenhanced rack counting invariants.
\end{abstract}

\textsc{Keywords:} Knots and links, Coxeter racks, finite racks,
knot and link invariants

\textsc{2000 MSC:} 57M25, 57M27, 17D99

\section{\large \textbf{Introduction}}

In 1982, Joyce introduced the term ``quandle'' for the algebraic
structure defined by translating the Reidemeister moves into
axioms for a binary operation corresponding to one strand crossing 
under another \cite{J}. In 1992, Fenn and Rourke generalized the 
category of quandles to the larger category of ``racks,'' whose 
axioms derive from framed isotopy moves \cite{FR}. The same
concepts appear under different names in the literature; quandles
are called ``distributive groupoids'' in \cite{M}, a special case
of quandle is called ``kei'' in \cite{T} and racks are described
as ``automorphic sets'' in \cite{B}.

One type of rack structure described in \cite{FR} is a \textit{Coxeter rack}, 
the subset of an $\mathbb{R}$-vector space on which a symmetric bilinear 
form is nonzero, with rack operation defined as a kind of reflection. The 
operator groups of such racks are related to Coxeter groups; see \cite{FR} 
for more.

Replacing $\mathbb{R}$ with $\mathbb{Z}_n$ yields finite Coxeter racks, 
which are suitable for use as target racks for counting invariants of 
knots and links as described in \cite{N2}. 
In \cite{NN} and \cite{NR}, quandle counting invariants are 
enhanced by making use of the module structure of the coloring 
quandles and biquandles, resulting in knot and link invariants which 
contain more information than the counting invariants alone. In this 
paper, we will apply the same idea to the counting invariants defined 
by finite Coxeter racks, obtaining a new family of enhanced rack counting 
invariants of knots and links.

The paper is organized as follows. In section \ref{rq} we review racks
and the rack counting invariant. In section \ref{crq} we review Coxeter 
racks and make a few observations. In section \ref{inv} we define the enhanced
Coxeter rack invariants and give examples demonstrating that the enhanced
invariants are strictly stronger than the unenhanced rack counting invariants. 
In section \ref{Q} we collect some questions for future research.

\section{\large \textbf{Racks and quandles}}\label{rq}

\begin{definition}
\textup{A \textit{rack} is a set $X$ with a binary operation 
$\tr:X\times X\to X$ such that
\begin{list}{}{}
\item[(i)]{for every pair $x,y\in X$ there is a unique $z\in X$ such that
$x=z\tr y$, and}
\item[(ii)]{for every triple $x,y,z\in X$ we have $(x\tr y)\tr z
=(x\tr z)\tr(y\tr z).$}
\end{list}}
\end{definition}

Rack axiom (i) says that every element $y\in X$ acts on $X$ via
a bijection $f_y:X\to X$, $f_y(x)=x\tr y$; the inverse defines a second
operation $x\tr^{-1} y = f^{-1}_y(x)$, and we have $(x\tr y)\tr^{-1} y =x$
and $(x\tr^{-1} y)\tr y =x$ for all $x,y\in X$. 
Rack axiom (ii) says that the operation $\tr$ is self-distributive. A rack
in which every element is idempotent, i.e., such that $x\tr x=x$ for all
$x\in X$, is a \textit{quandle}.

Rack structures abound in mathematics. Any algebraic structure which
acts on itself by automorphisms is a rack: define $f_y(x)=x\tr y$. Then
\[f_z(x\tr y) = f_z(x)\tr f_z(y) \longleftrightarrow
(x\tr y) \tr z = (x\tr z) \tr (y\tr z).\]
Standard examples of racks include:
\begin{list}{$\bullet$}{}
\item{A set $S$ with $x\tr y=\sigma(x)$ for some fixed bijection 
$\sigma:S\to S$ (\textit{constant action rack} or \textit{permutation rack})}
\item{A group $G$ with $x\tr y= y^{-n}xy^n$ (\textit{conjugation rack})}
\item{A group $G$ with $x\tr y = s(xy^{-1})y$ for a fixed automorphism
$s\in \mathrm{Aut}(G)$}
\item{A module over $\mathbb{Z}[t^{\pm 1},s]/s(s+t-1)$ with $x\tr y = tx + sy$
(\textit{$(t,s)$-rack})}
\end{list}

It is convenient to express a rack structure on a finite set
$\{x_1,x_2,\dots,x_n\}$ by encoding its
operation table as an $n\times n$ matrix $M$ whose $i,j$ entry is
$k$ where $x_k=x_i\tr x_j$. We call this matrix the \textit{rack matrix}
of $T$, denoted $M_T$.

\begin{example}
\textup{Let $R=\{1,2,3,4\}$ and $\sigma:R\to R$ be the permutation $(12)(34)$.
Then the rack matrix of $(R,\sigma)$ is}
\[M_{(12)(34)}=\left[\begin{array}{cccc}
2 & 2 & 2 & 2 \\
1 & 1 & 1 & 1 \\
4 & 4 & 4 & 4 \\
3 & 3 & 3 & 3
\end{array}\right].\]
\end{example}

For defining invariants of framed oriented knots and links, we need the 
\textit{fundamental rack}. By using the blackboard framing, we can consider link
diagrams as framed link diagrams with framing numbers given by the self-writhe
of each component; by choosing an order on the components, we can conveniently
express the writhe of an $n$-component link diagram as a vector 
$\mathbf{w}\in \mathbb{Z}^n$.

The idea is then to think of arcs in an oriented link diagram as generators
and the operation $\tr$ as ``crosses under from right to left'' when looking
in the positive direction of the overcrossing strand. The inverse
operation $\tr^{-1}$ then means ``crosses under from left to right.''
\[ \includegraphics{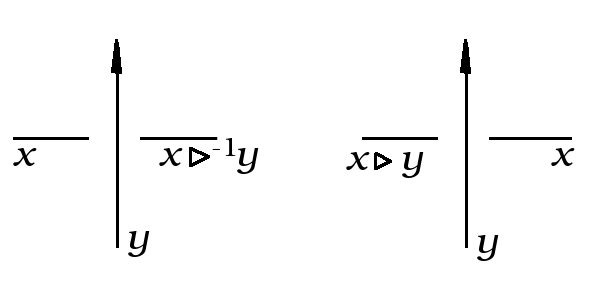}\] 
Indeed, the rack axioms are simply the Reidemeister moves required for framed 
isotopy interpreted in light of this operation; see \cite{FR} or \cite{N2} 
for more.

Given a link diagram $L$, we obtain a rack presentation with one 
generator for each arc and one relation at each crossing. That is, the 
\textit{fundamental rack} $FR(L)$
of the framed link specified by the diagram is the set of equivalence
classes of rack words under the equivalence relation generated by the crossing
relations together with the rack axioms. Note that changing the writhe of 
the diagram by Reidemeister I moves results in a generally different 
fundamental rack. Taking the quotient of $FR(L)$ for any framing of $L$ by 
setting $a\sim a\tr a$
for all $a\in FR(L)$ yields the \textit{knot quandle} of $L$, denoted $Q(L)$.

\begin{example} \textup{The pictured trefoil knot $3_1$ with writhe 3 has the 
listed fundamental rack presentation.}
\[\raisebox{-0.5in}{\includegraphics{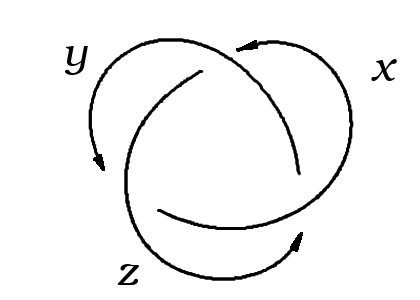}} \quad 
FR(D)= \langle x,\ y,\ z \ | \ x\tr y = z,\ y\tr z = x,\ z\tr x = y \rangle\]
\end{example}

In \cite{N2}, the quandle counting invariant $|\mathrm{Hom}(Q(K),T)|$ is 
extended to include non-quandle racks as coloring objects. For a finite rack
$T$, the \textit{rack rank} of $T$, denoted $N(T)$, is the exponent of the
permutation given by the diagonal of the rack matrix. If two ambient isotopic 
link diagrams $L$ and $L'$ have framing vectors with respect to an ordering 
on the components which are componentwise congruent mod $N(T)$, then there is
a bijection between the sets of rack colorings of $L$ and $L'$ by $T$. 
We might say that as far as $T$ cares, the writhes of $L$ live in 
$(\mathbb{Z}_{N(T)})^n$.
We are thus able to reduce the infinite set of fundamental racks of framings 
of $L$ to get a finitely computable ambient isotopy invariant of $L$.

\begin{definition}
\textup{Let $L$ be a link with $n$ components, $T$ a finite rack with rack 
rank $N(T)$ and $W=(\mathbb{Z}_{N(T)})^n$.
The \textit{polynomial rack counting invariant} of $L$ with respect to $T$ is}
\[PR(L,T)=\sum_{\mathbf{w}\in W} 
|\mathrm{Hom}(FR(D,\mathbf{w}),T)|\prod_{i=1}^nq_i^{w_i}\]
\textup{where $(D,w)$ is a diagram of $L$ with writhe vector $\mathbf{w}\in W$
and $FR(D,\mathbf{w})$ is the fundamental rack of $(D,\mathbf{w}).$}
\end{definition}

\begin{example}
\textup{The Hopf link $H$ has rack counting polynomial
$PR(H,T)=4+4q_1+4q_2+8q_1q_2$ with respect to the rack $T$ with rack matrix}
\[M_T=\left[\begin{array}{ccccc}
1 & 1 & 2 & 2 \\
2 & 2 & 1 & 1 \\
4 & 4 & 4 & 4 \\
3 & 3 & 3 & 3
\end{array}\right]\]
\textup{as the tables of colorings show.}

\[
\begin{array}{cccc}
\scalebox{0.9}{\includegraphics{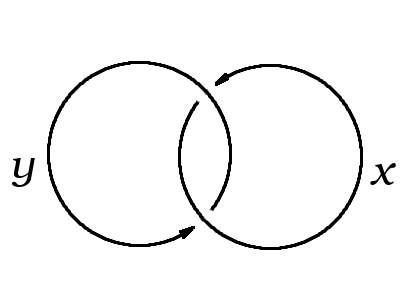}} &
\scalebox{0.9}{\includegraphics{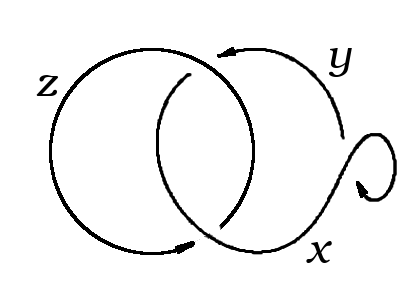}} &
\scalebox{0.9}{\includegraphics{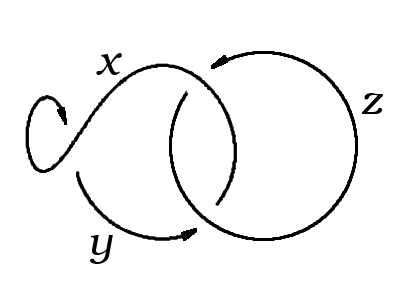}} &
\scalebox{0.9}{\includegraphics{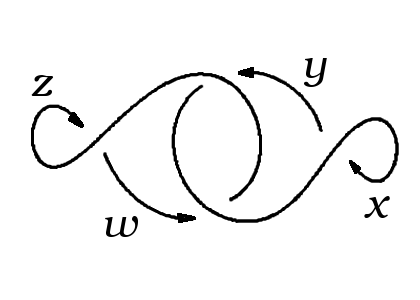}} \\
\begin{array}{|cc|} \hline
x & y \\ \hline
1 & 1 \\
1 & 2 \\
2 & 1 \\
2 & 2 \\ \hline
\end{array} &
\begin{array}{|ccc|} \hline
x & y & z\\ \hline
1 & 1 & 1 \\
1 & 1 & 2 \\
2 & 2 & 1 \\
2 & 2 & 2 \\ \hline
\end{array} &
\begin{array}{|ccc|} \hline
x & y & z \\ \hline
1 & 1 & 1 \\
1 & 1 & 2 \\
2 & 2 & 1 \\
2 & 2 & 2 \\ \hline
\end{array} &
\begin{array}{|cccc|} \hline
x & y & z & w \\ \hline
1 & 1 & 1 & 1 \\
1 & 1 & 2 & 2 \\
2 & 2 & 1 & 1 \\
2 & 2 & 2 & 2 \\ \hline
\end{array}\
\begin{array}{|cccc|} \hline
x & y & z & w \\ \hline
3 & 4 & 3 & 4 \\
3 & 4 & 4 & 3 \\
4 & 3 & 3 & 4 \\
4 & 3 & 4 & 3 \\
\hline
\end{array}
\end{array}\]

\textup{The two-component unlink $U_2$ has rack counting polynomial 
$PR(U_2,T)=16+8q_1+8q_2+4q_1q_2$ with respect to this rack as the reader
is invited to verify. Thus, the invariant detects the difference between
the Hopf link and the two-component unlink.}
\end{example}

\section{\large \textbf{Coxeter racks}}\label{crq}

In \cite{FR}, a \textit{Coxeter rack} is defined as the subset of
$\mathbb{R}^n$ on which a symmetric bilinear form $\langle,\rangle:
\mathbb{R}^n\times \mathbb{R}^n\to \mathbb{R}$ is nonzero, with rack operation
\[\mathbf{x}\tr \mathbf{y} = 
\mathbf{x}-\frac{2\langle\mathbf{x},\mathbf{y}\rangle}{\langle\mathbf{y},
\mathbf{y}\rangle}\mathbf{y}. \]
Multiplying the right-hand side by $-1$ defines 
a quandle structure, called a \textit{Coxeter quandle}.

We will study a family of slight generalizations of these Coxeter racks.
First, we replace $\mathbb{R}$ with an arbitrary commutative ring $R$,
choose a symmetric bilinear form $\langle ,\rangle:R^n\times R^n\to R$
and consider the subset of $R^n$ given by
\[T=\{\mathbf{x}\in R^n\ |\ \langle \mathbf{x},\mathbf{x}\rangle \in R^* \}\]
where $R^*$ is the set of units in $R$. Next, we note that replacing the $-1$ 
factor in the Coxeter quandle definition with any invertible scalar 
$\alpha\in R^*$ yields a valid rack structure, which gives us our generalized 
Coxeter rack definition. More formally, we have:

\begin{definition}
\textup{Let $R$ be a commutative ring, $V$ an $R$-module and 
$\langle,\rangle:V\times V\to R$ a symmetric bilinear form. Let }
\[T=\{\mathbf{x}\in V\ |\ \langle \mathbf{x},\mathbf{x}\rangle \in R^* \},\]
\textup{where $R^*$ is the set of units in $R$. Note that $T\subsetneq V$ since
$\mathbf{0}\not\in T$.
Define $\tr:T\times T\to T$ by}
\[\mathbf{x}\tr \mathbf{y} = \alpha\left(
\mathbf{x}-\frac{2\langle\mathbf{x},\mathbf{y}\rangle}{\langle\mathbf{y},
\mathbf{y}\rangle}\mathbf{y}\right)\]
\textup{where $\alpha\in R^*$. We call $(T,\tr)$ a \textit{generalized
Coxeter rack} and write $T=CR(V,\alpha, \langle,\rangle)$.}
\end{definition}

To verify that $(T,\tr)$ is a rack, we will find the following lemma useful:

\begin{lemma}
Let $(T,\tr)$ be a generalized Coxeter rack. Then 
$\langle\mathbf{x}\tr\mathbf{z},\mathbf{y}\tr\mathbf{z}\rangle = 
\alpha^2 \langle \mathbf{x},\mathbf{y}\rangle.$
\end{lemma}

\begin{proof}
\begin{eqnarray*}
\langle\mathbf{x}\tr\mathbf{z},\mathbf{y}\tr\mathbf{z}\rangle 
& = &
\left\langle\alpha\mathbf{x}
-\frac{2\alpha\langle\mathbf{x},\mathbf{z}\rangle}{\langle\mathbf{z},
\mathbf{z}\rangle}\mathbf{z},
\alpha\mathbf{y}
-\frac{2\alpha\langle\mathbf{y},\mathbf{z}\rangle}{\langle\mathbf{z},
\mathbf{z}\rangle}\mathbf{z}\right\rangle \\
& = &
\alpha^2\langle \mathbf{x},\mathbf{y}\rangle 
-2\alpha^2\frac{\langle\mathbf{y},\mathbf{z}\rangle}{\langle\mathbf{z},
\mathbf{z}\rangle}\langle\mathbf{x},\mathbf{z}\rangle
-2\alpha^2\frac{\langle\mathbf{x},\mathbf{z}\rangle}{\langle\mathbf{z},
\mathbf{z}\rangle}\langle\mathbf{z},\mathbf{y}\rangle
+4\alpha^2\frac{\langle\mathbf{x},\mathbf{z}\rangle}{\langle\mathbf{z},
\mathbf{z}\rangle}
\frac{\langle\mathbf{y},\mathbf{z} \rangle}{\langle\mathbf{z},\mathbf{z}\rangle}\langle \mathbf{z},\mathbf{z}\rangle \\
& = &
\alpha^2\langle \mathbf{x},\mathbf{y}\rangle. \\
\end{eqnarray*}
\end{proof}

\begin{proposition}
$T$ is a rack under $\tr$.
\end{proposition}

\begin{proof}

To see that $\tr$ is right-invertible, note that the operation 
$\tr^{-1}:M\times M\to M$ defined by 
\[\mathbf{x}\tr^{-1} \mathbf{y} = \alpha^{-1}\left(
\mathbf{x}-\frac{2\langle\mathbf{x},\mathbf{y}\rangle}{\langle\mathbf{y},
\mathbf{y}\rangle}\mathbf{y}\right)\]
satisfies $(\mathbf{x}\tr \mathbf{y})\tr^{-1} \mathbf{y}=\mathbf{x}$:
\begin{eqnarray*}
(\mathbf{x}\tr \mathbf{y})\tr \mathbf{y} & = & 
\alpha^{-1}\left(\mathbf{x}\tr\mathbf{y} 
-\frac{2}{\langle\mathbf{y},\mathbf{y}\rangle}
\langle\mathbf{x}\tr\mathbf{y},\mathbf{y}\rangle\mathbf{y}\right) \\
& = & 
\alpha^{-1}\left(\alpha\left(
\mathbf{x}-\frac{2\langle\mathbf{x},\mathbf{y}\rangle}{\langle\mathbf{y},
\mathbf{y}\rangle}\mathbf{y}\right)
-\frac{2}{\langle\mathbf{y},\mathbf{y}\rangle}
\left\langle\alpha\left(
\mathbf{x}-\frac{2\langle\mathbf{x},\mathbf{y}\rangle}{\langle\mathbf{y},
\mathbf{y}\rangle}\mathbf{y}\right),\mathbf{y}\right\rangle\mathbf{y}\right) \\
& = &
\mathbf{x} -\frac{2\langle\mathbf{x},\mathbf{y}\rangle}{\langle\mathbf{y},
\mathbf{y}\rangle}\mathbf{y}
-\frac{2\langle\mathbf{x},\mathbf{y}\rangle}{\langle\mathbf{y},
\mathbf{y}\rangle}\mathbf{y}
+\frac{4\langle\mathbf{x},\mathbf{y}\rangle \langle\mathbf{y},
\mathbf{y}\rangle}{\langle\mathbf{y},\mathbf{y}\rangle^2}\mathbf{y} \\
& = & \mathbf{x}.
\end{eqnarray*}
In particular, note that if $\alpha=\alpha^{-1}$ then $\tr=\tr^{-1}$ and 
$(M,\tr)$ is an \textit{involutory} rack.

Finally, we check that $\tr$ is self-distributive:
\begin{eqnarray*}
(\mathbf{x}\tr \mathbf{y}) \tr \mathbf{z} & = & \alpha(\mathbf{x}\tr \mathbf{y}) 
-2\alpha \frac{\langle \mathbf{x}\tr \mathbf{y},\mathbf{z}\rangle}{\langle \mathbf{z},\mathbf{z} \rangle}\mathbf{z} \\
& = & \alpha^2 \mathbf{x}-2\alpha^2\frac{\langle \mathbf{x},\mathbf{y} \rangle}{\langle \mathbf{y},\mathbf{y}\rangle}  \mathbf{y} 
+\left(-2\alpha^2 \frac{\langle \mathbf{x},\mathbf{z}\rangle}{\langle \mathbf{z},\mathbf{z} \rangle}
+4\alpha^2\frac{\langle \mathbf{x},\mathbf{y}\rangle\langle \mathbf{y},\mathbf{z} \rangle}{\langle \mathbf{y},\mathbf{y}\rangle\langle \mathbf{z},\mathbf{z}\rangle}\right)\mathbf{z} \\
\end{eqnarray*}
On the other hand,
\begin{eqnarray*}
(\mathbf{x}\tr \mathbf{z})\tr (\mathbf{y}\tr \mathbf{z}) & = &
\alpha(\mathbf{x}\tr \mathbf{z}) -2\alpha\frac{\langle \mathbf{x}\tr \mathbf{z}, \mathbf{y}\tr \mathbf{z}\rangle}{\langle \mathbf{y}\tr \mathbf{z}, \mathbf{y}\tr \mathbf{z} \rangle}(\mathbf{y}\tr \mathbf{z}) \\
& = &
\alpha(\mathbf{x}\tr \mathbf{z}) -2\alpha\frac{\langle \mathbf{x}, \mathbf{y}\rangle}{\langle \mathbf{y}, \mathbf{y} \rangle}(\mathbf{y}\tr \mathbf{z}) \\
& = &
\alpha^2\mathbf{x}-2\alpha^2\frac{\langle \mathbf{x},\mathbf{z} \rangle}{\langle \mathbf{z},\mathbf{z}\rangle} \mathbf{z} 
-2\alpha\frac{\langle \mathbf{x}, \mathbf{y}\rangle}{\langle \mathbf{y}, \mathbf{y} \rangle}
\left(\alpha \mathbf{y}-2\alpha\frac{\langle \mathbf{y},\mathbf{z} \rangle}{\langle \mathbf{z},\mathbf{z}\rangle} \mathbf{z}\right) \\
& = & \alpha^2 \mathbf{x}-2\alpha^2\frac{\langle \mathbf{x},\mathbf{y} \rangle}{\langle \mathbf{y},\mathbf{y}\rangle}  \mathbf{y} 
+\left(-2\alpha^2 \frac{\langle \mathbf{x},\mathbf{z}\rangle}{\langle \mathbf{z},\mathbf{z} \rangle}
+4\alpha^2\frac{\langle \mathbf{x},\mathbf{y}\rangle\langle \mathbf{y},\mathbf{z} \rangle}{\langle \mathbf{y},\mathbf{y} \rangle\langle \mathbf{z},\mathbf{z}\rangle}\right)\mathbf{z} \\
\end{eqnarray*}
as required.
\end{proof}

\begin{remark}\textup{
In the classical case where $V=\mathbb{R}^n$, $\langle,\rangle$ is
the dot product and $\alpha=-1$, the Coxeter quandle operation is the 
result of reflecting $\mathbf{x}$ through $\mathbf{y}$ in the plane 
spanned by $\mathbf{x}$ and $\mathbf{y}$:}
\[\includegraphics{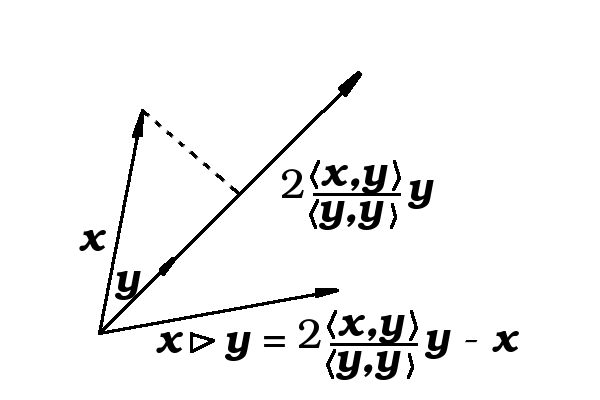}\]
\end{remark}

For the purpose of defining counting invariants, we need finite racks.
Thus, we will consider the case $V=(\mathbb{Z}_n)^m$ of generalized Coxeter
racks which are subsets of free modules over the integers modulo $n$. 
The input data required to construct such a rack consists of two
integer parameters $n$ and $m$ (which determine $V=(\mathbb{Z}_n)^m$),
a scalar $\alpha\in \mathbb{Z}_n$ coprime to $n$, and
a symmetric matrix $A\in M_m(\mathbb{Z}_n)$ which defines a symmetric bilinear
form $\langle,\rangle:(\mathbb{Z}_n)^m\times (\mathbb{Z}_n)^m\to \mathbb{Z}_n$
by \[\langle\mathbf{x},\mathbf{y}\rangle=\mathbf{x}A\mathbf{y}^t\]
where $\mathbf{x}$ is a row vector and $\mathbf{y}^t$ is a column vector.

\begin{example}
\textup{For a simple example, let us take $V=(\mathbb{Z}_3)^2$ with
$\alpha=1$ and $A=\left[\begin{array}{cc}1 & 2 \\ 2 & 0 \end{array}\right]$.
A straightforward computation then shows that 
$T=CR(V,\alpha,A)=\{ x_1=(1,0),\ x_2=(1,1),\ x_3=(2,0),\ x_4=(2,2)\}$ and 
T has rack matrix}
\[M_T=\left[\begin{array}{cccc}
3 & 1 & 3 & 1 \\
2 & 4 & 2 & 4 \\
1 & 3 & 1 & 3 \\
4 & 2 & 4 & 2 \\
\end{array}\right].\]
\end{example}

We end this section with a few brief observations about finite
Coxeter racks.

\begin{proposition}
If $R=\mathbb{Z}_2$, then every generalized Coxeter rack 
over $R$ has trivial rack operation.
\end{proposition}

\begin{proof}
\[\mathbf{x}\tr \mathbf{y} = 1\left(
\mathbf{x}-\frac{2\langle\mathbf{x},\mathbf{y}\rangle}{\langle\mathbf{y},
\mathbf{y}\rangle}\mathbf{y}\right) = \mathbf{x} -0 =\mathbf{x}.\]
\end{proof}

\begin{corollary}
Finite generalized Coxeter racks over $R=\mathbb{Z}_2$ are 
classified by the cardinality of the subset of $V$ on which $\langle,\rangle$ 
is nonzero.
\end{corollary}

\begin{proposition}
If $R$ has characteristic 2, then 
$CR(R^m,\alpha,A)$ is a constant action rack with permutation given by 
multiplication by $\alpha.$
\end{proposition}

\begin{proof}
If the characteristic of $R$ is 2, we have
\[\mathbf{x}\tr\mathbf{y} = \alpha\mathbf{x}+0 =\alpha\mathbf{x}\]
for all $\mathbf{x}\in T$.
\end{proof}

\begin{proposition} \label{creq}
Let $T=CR(R^m,\alpha, A)$ and $T'=CR(R^m,\alpha,\beta A)$ for an
invertible scalar $\beta\in R^*$. Then $T=T'$.
\end{proposition}

\begin{proof}
First, $\mathbf{x}A\mathbf{x}^t\in R^*$
if and only if $\mathbf{x}\beta A\mathbf{x}^t\in R^*$ so setwise
$T=T'$. Moreover, 
\[ \mathbf{x}\tr_T \mathbf{y}  = 
\alpha\left( \mathbf{x}
-\frac{2\mathbf{x}A\mathbf{y}^t}{\mathbf{y}A\mathbf{y}^t}\mathbf{y}\right)
 = 
\alpha\left(\mathbf{x}-\frac{2\mathbf{x}\beta A\mathbf{y}^t}{
\mathbf{y}\beta A\mathbf{y}^t}\mathbf{y}\right)
 = (\mathbf{x}\tr_{T'} \mathbf{y})
\]
so the rack structures are the same.
\end{proof}

\begin{remark}
\textup{Changing the scalar $\alpha$ in general does result in different
rack structures: $\alpha=-1$ yields a quandle while $\alpha=1$ yields
a non-quandle rack, as we know. As another example, we note that the 
Coxeter racks 
$CR\left((\mathbb{Z}_5)^2,3,\left[\begin{array}{cc}1 & 2 \\ 2 & 0 
\end{array}\right]\right)$ and 
$CR\left((\mathbb{Z}_5)^2,1,\left[\begin{array}{cc}1 & 2 \\ 2 & 0 
\end{array}\right]\right)$ respectively have rack polynomials $16$ and 
$16s^4t^4$ (see \cite{N}), and hence cannot be isomorphic.}

\end{remark}

\section{\large \textbf{Coxeter enhanced rack counting invariants}}\label{inv}

Let $R=\mathbb{Z}_n,$ $V=R^m$ and $T=CR(V,\alpha,A)$. Let $L$ be a link with
$c$ components and let $W=(\mathbb{Z}_{N(T)})^c$. Consider a coloring 
$f\in \mathrm{Hom}(FR(D,\mathbf{w}),T)$ of an oriented framed link diagram
$(D,\mathbf{w})$. For each homomorphism $f$, the image subrack 
$\mathrm{Im}(f)$ and hence the submodule $\mathrm{Span}(\mathrm{Im}(f))
\subseteq V$ it spans are invariant under Reidemeister moves, so we can form an
enhanced version of the counting invariant incorporating this
extra information. Formally, we have

\begin{definition}
\textup{Let $L$ be a link with $n$ components and $T=CR(R,T,\alpha,A)$
a finite generalized Coxeter rack. Then the \textit{Coxeter enhanced
rack counting invariant} is}
\[cp(L,T) = \sum_{\mathbf{w}\in W} 
\left(\sum_{f\in\mathrm{Hom}(FR(D,\mathbf{w}),T)}
\prod_{i=1}^nq_i^{w_i}s^{|\mathrm{Span}(\mathrm{Im}(f))|}t^{|\mathrm{Im}(f)|}\right).\]
\end{definition}

We note from the definition that specializing $s=t=1$ yields the rack counting
polynomial. On the other hand, the Coxeter enhanced invariants are stronger 
than the corresponding unenhanced rack counting invariants, as the following 
example shows.

\begin{example}
\textup{Let $L_1$ be the $(4,2)$ torus link and $L_2$ the three-component link
illustrated below. Let $T$ be the Coxeter quandle $T=CR\left((\mathbb{Z}_3)^2,2,
\left[\begin{array}{cc} 1 & 2 \\ 2 & 0\end{array}\right]\right)$. Then $T$ has 
quandle matrix}
\[M_T=\left[\begin{array}{cccc}
1 & 3 & 1 & 3 \\
4 & 2 & 4 & 2 \\
3 & 1 & 3 & 1 \\
2 & 4 & 2 & 4 \\
\end{array}\right]\]
\textup{where we have $x_1=(1,0),\ x_2=(1,1),\ x_3=(2,0)$ and $x_4=(2,2)$.
An easy computation then shows that while both $L_1$ and $L_2$ have
quandle counting invariant $|\mathrm{Hom}(Q(L_i),T)|=16$, the
Coxeter enhanced invariants tell the links apart, with 
$cp(L_1,T)=4s^3t+4s^3t^2+8s^9t^4$ while 
$cp(L_2,T)=4s^3t+12s^3t^2$.}
\[\begin{array}{cc}
\includegraphics{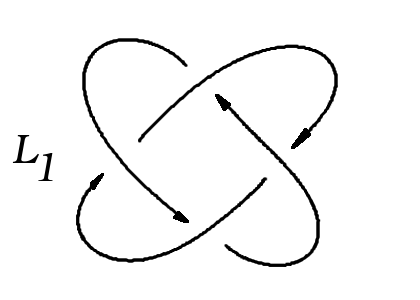} & \includegraphics{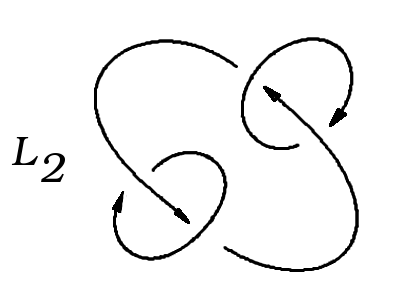} \\
4s^3t+4s^3t^2+8s^9t^4 & 4s^3t+12s^3t^2.
\end{array}
\]
\end{example}

In light of proposition \ref{creq}, we have
\begin{proposition}
For any link $L$, $cp(L,CR(V,\alpha,A))=cp(L,CR(V,\alpha,\beta A))$
for any invertible scalar $\beta\in R^*$.
\end{proposition}

\begin{example}
\textup{The Coxeter enhanced rack counting invariants give us convenient
information about the subracks of $T=CR(V,\alpha,A)$ with surjective
homomorphisms from the various fundamental racks of the framings of $L$.
The trefoil knot $3_1$ has Coxeter enhanced rack invariant 
$cp=6s^3t^2+12s^9t^6$ with respect to the Coxeter rack 
$T=\left((\mathbb{Z}_3)^2,1,
\left[\begin{array}{cc} 1 & 1 \\ 1 & 1 \end{array}\right]\right)$.
We can immediately note several things: there are no homomorphisms
from the fundamental rack of trefoil with odd framing into $T$ since the
coefficient of $q^1$ is zero; there are six surjective homomorphisms
from the fundamental rack of any evenly-framed trefoil onto subracks of
$T$ with two elements spanning 1-dimensional subspaces of 
$(\mathbb{Z}_3)^2$, and there are twelve surjective homomorphisms
from the fundamental rack of an evenly framed trefoil onto subracks
of $T$ with six elements spanning the whole space $(\mathbb{Z}_3)^2$
(as indeed six vectors in a two-dimensional subspace must).}

\end{example}

Finally, we note that to define the Coxeter enhanced rack counting 
invariants, we need to know the vector space or $R$-module structure
of the Coxeter rack in question, not just the rack structure. It
seems possible \textit{a priori} that the same Coxeter rack might
embed in different modules or vector spaces or in different ways in 
the same module or vector space (indeed, see the next section), and in
such a case we expect the resulting invariants to be different, although 
of course they specialize to the same rack counting invariant.

\section{\large \textbf{Questions}}\label{Q}

In this section we collect some questions for future research.

Given a finite Coxeter rack $T$, to what degree is it possible to
recover the module structure $(R,V,\alpha,A)$ from which 
$T$ arises? We already know that changing $A$ by an invertible scalar 
multiple yields an isomorphic (indeed, an \textit{identical}) rack; an
easy computation shows that the Coxeter racks on $(\mathbb{Z}_3)^2$ with
$\alpha=1$ and bilinear forms defined by 
$A=\left[\begin{array}{cc} 1 & 2 \\ 2 & 0 \end{array}\right]$ and
$B=\left[\begin{array}{cc} 0 & 2 \\ 2 & 0 \end{array}\right]$ are not
identical but are nonetheless isomorphic.

In \cite{FR} we also find \textit{Hermitian form} racks, which are 
defined on the subset of a vector space on which 
$\langle\mathbf{x},\mathbf{x}\rangle\in R^*$ for a form 
$\langle,\rangle:V\times V\to R$ which is linear in the first variable and
conjugate linear in the second for a conjugation on $R$ (that is,
a ring automorphism $f:R\to R$ such that $f(f(r))=r$ for all $r\in R$). Our
generalized Coxeter racks are just the Hermitian racks on 
$\mathbb{Z}_n$-modules with respect to the identity conjugation, 
the only conjugation on rings whose additive groups are cyclic.
We have not considered Hermitian racks on modules over finite rings
other than $\mathbb{Z}_n$, but the same construction given in section
\ref{inv} should give a \textit{Hermitian enhanced rack counting invariant}
for any finite $R$-module with a nontrivial conjugation, e.g. vector spaces 
over the Galois field of four elements with the Frobenius automorphism.

What kinds of non-rack birack structures can be defined on finite vector
spaces and modules with Coxeter/Hermitian type operations? Such structures 
should also give rise to enhanced counting invariants.

\bigskip

Our \texttt{python} code for computing Coxeter rack matrices and Coxeter
enhanced rack counting invariants is available at \texttt{www.esotericka.org.}

\bigskip

\textsc{Department of Mathematics, Claremont McKenna College,
 850 Colubmia Ave., Claremont, CA 91711}

\noindent
\textit{Email address: }\texttt{knots@esotericka.org}

\medskip

\textsc{Department of Mathematics, Pomona College,
 610 N. College Ave., Claremont, CA 91711}

\noindent
\textit{Email address: }\texttt{Ryan.Wieghard@pomona.edu}


\begin{thebibliography}{0}

\bibitem{B}{E. Brieskorn, Automorphic sets and braids and singularities. 
\textit{Contemp. Math.} \textbf{78} (1988) 45-115.} 

\bibitem{J}{D. Joyce.
 A classifying invariant of knots, the knot quandle.
 \textit{J. Pure Appl. Algebra}  \textbf{23}  (1982)  37-65.}

\bibitem{FR}{R. Fenn and C. Rourke.
 Racks and links in codimension two.
 \textit{J. Knot Theory Ramifications}  \textbf{1}  (1992), 343-406.}

\bibitem{M}{S. V. Matveev.
Distributive groupoids in knot theory. 
\textit{Math. USSR, Sb.} \textbf{47} (1984) 73-83.}

\bibitem{NN}{E.A. Navas and S. Nelson. On symplectic quandles. 
To appear in \textit{Osaka J. Math.}, arXiv:math/0703727}

\bibitem{N}{S. Nelson. A polynomial invariant of finite racks.
\textit{J. Alg. Appl.} \textbf{7} (2008) 263-273.}

\bibitem{N2}{S. Nelson. Link invariants from finite racks.
arXiv:0808.0029}

\bibitem{NR}{S. Nelson and J.L. Rische. On bilinear biquandles.
\textit{Colloq. Math.} \textbf{112} (2008) 279-289.}

\bibitem{T}{M. Takasaki. Abstraction of symmetric transformation 
(in Japanese). \textit{Tohoku Math J.} \textbf{49} (1943) 145-207.}


\end{thebibliography}
\end{document}